\documentclass{amsart}

\usepackage[colorlinks=true, urlcolor=blue,bookmarks=true,bookmarksopen=true,citecolor=blue,hypertex]{hyperref}
\usepackage{graphicx}
\usepackage{enumerate}

\usepackage{amscd}
\usepackage{amssymb}
\usepackage{amsxtra}
\usepackage{amsmath}
\usepackage{enumerate}
\usepackage{mathrsfs}

\usepackage{color}

\usepackage{amsmath,amsfonts,amssymb}
\usepackage{verbatim}
\usepackage[hmargin = 4cm,vmargin = 2.5cm]{geometry}
\usepackage{epsfig}
\usepackage{amsthm}
\usepackage{bbold}

\usepackage[T1]{fontenc}

\theoremstyle{plain}

\newtheorem{thm}{Theorem}

\newtheorem{cor}[thm]{Corollary}

\theoremstyle{definition}

\theoremstyle{remark}
\newtheorem{rem}[thm]{Remark}

\theoremstyle{plain}

\def\NN{\mathbb{N}}
\def\ZZ{\mathbb{Z}}

\def\RR{\mathbb{R}}

\def\dd{\mathrm{d}}
\def\Ham{Ham}
\def\Symp{Symp}
\def\Cal{Cal}

\def\Disk{D}

\def\Area{\operatorname{Area}}

\begin{document}

\pagestyle{headings}

\bibliographystyle{alphanum}

\title{Non-autonomous curves on surfaces}

\thanks{The author was supported by the Azrieli Fellowship.}

\date{\today} 
\author{Michael Khanevsky}
\address{Michael Khanevsky, Mathematics Department,
Technion - Israel Institute of Technology
Haifa, 32000,
Israel}
\email{khanev@math.technion.ac.il}

\begin{abstract}
Consider a symplectic surface $\Sigma$ with two properly embedded Hamiltonian isotopic curves $L$ and $L'$. Suppose $g \in \Ham (\Sigma)$ is a Hamiltonian diffeomorphism
which sends $L$ to $L'$. Which dynamical properties of $g$ can be detected by the pair $(L, L')$? We present two scenarios where one can deduce that $g$ is `chaotic': non-autonomous
or even of positive entropy.
\end{abstract}

\maketitle

\section{Introduction and results}

Given a Hamiltonian diffeomorphism $g$ it is extremely difficult to analyze it. It can be seen already at the stage of extracting numerical information: most of useful invariants
(e.g. entropy, spectral data related to periodic points, etc.) are not easy to compute in the general case. Instead of attacking $g$ itself one may consider its action on spaces
that are easier to understand. We restrict our attention to Hamiltonian diffeomorphisms on surfaces and their action on curves (Lagrangian submanifolds).

Clearly, given a pair of properly embedded curves $L$ and $L' = g(L)$ on a surface $\Sigma$ it is easy to extract certain numerical data: the symplectic area of connected components
of $\Sigma \setminus (L \cup L')$ or combinatorial data associated to their partition of $\Sigma$. In fact, in generic situation, this gives a complete set of invariants: one can reconstruct
the pair $(L, L')$ up to a diagonal action by a symplectomorphism. That is, up to a symplectic change of coordinates.
The main question is to what extent the behavior of $g$ can be detected by looking at $L$ and $L'$ rather than at $g$ itself. For example, 
in ~\cite{RSS:combi-floer-homology} the authors show how data described above can be used to compute the Lagrangian Floer homology of $(L, L')$ which, in turn, has well-established relation to
the Floer-theoretic data of $g$. In this article we describe scenarios where $(L, L')$ provide evidence that $g$ is `chaotic' -- has positive topological entropy or at least is non-autonomous. 

We prove the following.

\begin{thm}\label{t:entropy}
  Suppose $\Sigma$ is a compact connected symplectic surface, possibly with boundary and punctures, and $L$ is an essential simple closed curve in $\Sigma$. Pick $h_0 > 0$ to be a threshold on entropy.
  Then there exists a curve $L_{h_0}$ Hamiltonian isotopic to $L$ which satisfies the following. For every $g \in \Ham (\Sigma)$ such that $g (L) = L_{h_0}$, the topological entropy
  $h (g) > h_0$.
\end{thm}

\begin{cor} 
  We define the topological entropy of a pair of essential Hamiltonian isotopic curves:
  \[
	  h (L', L) = \inf \{h (g) \; \big| \; g \in \Ham (\Sigma) \; \text{s.t.} \; g (L) = L'\}.
  \]
The theorem shows that this invariant is unbounded, in particular, not identically zero.
\end{cor}

The proof uses quasimorphisms on $\Ham (\Sigma)$ constructed by Brandenbursky and Marcin\-kowski ~\cite{Bra-Ma:ent-quas} that are Lipschitz with respect to the topological entropy.
We show that they descend to invariants of pairs of essential curves. These invariants are ill-defined in the sense that they can be computed up to a bounded ambiguity 
(the defect of the quasimorphism) but that is sufficient when one tries to analyze behavior on a large scale.
Using the same tools, one can show that the entropy metric on $\Ham(\Sigma)$ (word metric with respect to the generating set of entropy-zero Hamiltonians) or the autonomous metric
(word metric with respect to autonomous Hamiltonians) are not bounded in the orbit  $\{ g (L) \; \big| \; g \in \Ham (\Sigma)\}$ of an essential curve $L$. 

It would be interesting to obtain similar results on surfaces that do not admit essential curves (e.g. sphere, disk, annulus). ~\cite{Bra-Ma:ent-quas} provides
a large family of entropy-Lipschitz quasimorphisms, which, however, do not descend to the space curves. In the case of an annulus we use Calabi quasimorphisms
constructed by Entov and Polterovich ~\cite{En-Po:calqm} to show a somewhat weaker statement:

\begin{thm} \label{t:annulus}
  Let $\Sigma = S^1 \times [0, 1]$ be an annulus equipped with the standard symplectic form 
  and $L = \{0\} \times [0, 1] \cup \{\frac{1}{2}\} \times [0, 1]$ be a union of two chords (we use the convention $S^1 = \RR / \ZZ$). 
  There exists an $L'$ Hamiltonian isotopic to $L$ such that all $g \in \Ham (\Sigma)$ with $g (L) = L'$ are not autonomous.
  
  In addition, given $R > 0$, one may pick $L'$ such that the distance between $\{ g \; \big| \; g (L) = L' \}$ to the set of autonomous Hamiltonian diffeomorphisms 
  is at least $R$ in the Hofer metric.
\end{thm}

In this paper Hamiltonian diffeomorphisms are assumed to have compact support in the interior of $\Sigma$.
We need to stress that by \emph{autonomous} Hamiltonians we understand those generated by an autonomous flow supported in the interior of $\Sigma$. 
The proposed proof of the theorem fails if one extends the definition to flows that rotate the boundary of the annulus.

While this theorem states a weaker result than Theorem~\ref{t:entropy} from a dynamical point of view, it adds an important geometric perspective. 
In the proof we present an explicit construction of such $L'$.

Like before, the quasimorphisms descend to [ill-defined] invariants of pairs $(g L, L)$. In the case when $g$ is autonomous they provide information on the Reeb graph of a Hamiltonian function generating $g$.
This data can be compared with that extracted from the curves directly (e.g. estimates on rotation numbers of points in $\Sigma$ under $g$). 
In our example this will result in a contradiction which means that $g$ cannot be autonomous.

This example shows that the set $\{ (g L, L) \; \big| \; g \in \Ham (\Sigma)\}$ has diameter greater or equal to two in the autonomous metric. In fact, 
we construct $L'$ by deforming $L$ by two autonomous Hamiltonians. It would be interesting to find an example where the distance is at least three.

\medskip

\emph{Acknowledgements:}
The author wishes to thank M. Brandenbursky, M. Entov and L. Polterovich for their comments on these results. We are also grateful to 
the referee for remarks regarding organization of this text.

\section{Definitions}\label{S:def}
 
Let $G$ be a group. A function $r : G \to \RR$ is called a \emph{quasimorphism} if there exists 
a constant $D$ (called the \emph{defect} of $r$) such that $|r(fg) - r(f) - r(g)| < D$ for all $f, g \in G$. The quasimorphism $r$ 
is called \emph{homogeneous} if it satisfies $r(g^m) = mr(g)$ for all $g \in G$ and $m \in \ZZ$.
Any homogeneous quasimorphism satisfies $r(fg) = r(f) + r(g)$ for commuting elements $f, g$. Every quasimorphism
is equivalent (up to a bounded deformation) to a unique homogeneous one ~\cite{Ca:scl}. 

\medskip

Let $L$ be a curve in a symplectic surface $\Sigma$. $\Ham(\Sigma)$ denotes the group of Hamiltonian diffeomorphisms with compact support in the interior of $\Sigma$.
$$S = \{ g \in \Ham (\Sigma) \; | \; g (L) = L\}$$ is the stabilizer of $L$.
The orbit $O_L = \{ g (L) \; \big| \; g \in \Ham (\Sigma)\}$ can be identified with the set $\Ham (\Disk) / S$ of left cosets of $S$.

Let $r : \Ham (\Sigma) \to \RR$ be a quasimorphism which vanishes on $S$. Using the identification $O_L \simeq \Ham (\Disk) / S$ one shows that for all 
$g, h \in \Ham (\Sigma)$ such that $[g] = [h] \in O_L$, $g$ differs from $h$ by an element of $S$, hence $|r (g) - r (h)| \leq D$ ($D$ is the defect of $r$).
Consequently, $r$ induces an ill-defined function  $r_L : O_L \to \RR$. It can be treated either as a set-valued function whose values have bounded distribution or as a function
which is defined up to ambiguity $D$. Another option is to pick a representative in each coset. We will use the first alternative. In this case notation $r_L (L') > h$
means that all elements of the set $r_L (L')$ are greater than $h$.

A metric on $\Ham (\Sigma)$ induces a pseudo-metric on the orbit $O(L)$ by
\[
  d (L_1, L_2) = \inf \{\|g\| \, \big| \, g \in \Ham (\Sigma) \; \text{s.t.} \; g (L_1) = L_2\}.
\]
If the metric on $\Ham$ is discrete (like in the case of various word metrics), it induces a genuine metric on $O_L$.

\medskip

The set of autonomous Hamiltonian diffeomorphisms and the set of entropy-zero ones both generate $\Ham (\Sigma)$. The word metrics 
on $\Ham(\Sigma$) with respect to these generating sets are called the autonomous and the entropy metric.
Since any homogeneous quasimorphism which is Lipschitz with respect to the entropy vanishes on all autonomous and all entropy-zero
Hamiltonians, such quasimorphisms are Lipschitz also in the autonomous and the entropy metric.

\section{Unbounded entropy}

A simple closed curve $L \subset \Sigma$ is called \emph{essential} if it is not contractible, not isotopic to a boundary curve and cannot be contracted to a puncture.

We prove Theorem~\ref{t:entropy}. ~\cite{Bra-Ma:ent-quas} constructs an infinite-dimensional family of homogeneous quasimorphisms $r: \Ham (\Sigma) \to \RR$ which are Lipschitz with respect to the 
topological entropy:
\[
  |r (g)| \leq C_r h (g).
\]
Given an essential curve $L \subset \Sigma$, we show below that these quasimorphisms vanish on the stabilizer of $L$. Therefore, given $g, f \in \Ham (\Sigma)$ such that
$g L = f L$ it holds $$r (f) \geq r (g) - D_r.$$ Any non-trivial homogeneous quasimorphism $r$ is unbounded. Given $h_0 > 0$ pick $g$ with $r(g) > C_r h_0 + D_r$ and
put $L_{h_0} = g (L)$. It follows that all $f$ with $f L = g L = L_{h_0}$ satisfy $h (f) > h_0$.

\medskip

In other words, $r$ induces an ill-defined invariant $r_L : O_L \to \RR$ as explained in the previous section. $r_L$ is unbounded and $|r_L (L')|-D_r$ 
induces a lower bound for the entropy $h (L', L)$. At the same time the quasimorphisms from ~\cite{Bra-Ma:ent-quas} are Lipschitz also with respect to the autonomous and the entropy metrics.
This implies that $|r_L (\cdot)|-D_r$ provides lower bounds also for the induced metrics on $O_L$, hence $O_L$ has infinite diameter. Using the fact that the family of quasimorphisms $r$ is `large',
one can use standard arguments to deduce that $O_L$ admits quasi-isometric embeddings of `large' subsets (e.g. $\ZZ^N$ for all $N > 0$).

\begin{rem}
  The argument below works verbatim if one replaces $\Ham (\Sigma)$ with $\Symp (\Sigma)$ or $\Symp_0 (\Sigma)$, leading to the same results.
  Similarly, one may impose a lower bound in the autonomous or the entropy norm instead of the lower bound on entropy.
\end{rem}

\medskip

We briefly describe the construction of quasimorphisms in ~\cite{Bra-Ma:ent-quas},
while the reader is invited to consider the article for more detailed definitions and proofs.
Bestvina and Fujiwara constructed a family of quasimorphisms $\psi: MCG (\Sigma_n) \to \RR$ 
where $MCG (\Sigma_n)$ is the mapping class group of $n$-times punctured $\Sigma$ (~\cite{Be-Fu:bounded-cohomology}).
As every quasimorphism can be homogenized, we may assume that $\psi$ is homogeneous. 
Pick $n$ distinct points $\mathbf{z} = (z_1, \ldots, z_n)$ in the interior of $\Sigma$. Given an $n$-tuple $\mathbf{x} = (x_1, \ldots, x_n)$ in the configuration space $X_n (\Sigma)$, 
push each $z_j$ to $x_j$ by an isotopy supported near a short geodesic path, compose with $g$ and finally push
each $g (x_j)$ back to $z_j$ along a short geodesic path. 
One has to exclude certain problematic configurations (for example, those where a marked point lies on a geodesic path
or those featuring multiple length minimizing geodesics connecting points of interest). 
For the remaining $\mathbf{x} \in X_n (\Sigma)$ this construction results in a diffeomorphism $g_\mathbf{x}$ of $\Sigma$ which fixes the $n$-tuple of marked points
$\mathbf{z}$ and determines an element $\gamma (g, \mathbf{x}) = [g_\mathbf{x}] \in MCG (\Sigma_n)$.
We remark that though $g_\mathbf{x}$ is isotopic to the identity (hence it represents the identity element of $MCG (\Sigma)$), 
it may induce a highly nontrivial action on the topology of the punctured surface.

The quasimorphism $r$ is defined by
\[
	r (g) = \lim_{p \to \infty} \frac{1}{p} \int_{X_n} \psi (\gamma (g^p, \mathbf{x})) \dd \text{vol}.
\]
The set of problematic configurations $\mathbf{x}$ where the construction of $\gamma (g, \mathbf{x})$ fails is a null set 
hence can be ignored under the integration.
For an appropriate choice of $\psi: MCG (\Sigma_n) \to \RR$, $r$ is Lipschitz with respect to the entropy $h(g)$.

\medskip


Now pick an essential curve $L \in \Sigma$. We show that $r$ vanishes on the stabilizer of $L$. 

Suppose first that $\Sigma \setminus L$ is connected. 
Let $g$ be a Hamiltonian that fixes $L$. We pick the marked points $\mathbf{z}$ away from $L$ and restrict attention to the configuration 
space $X_n (\Sigma \setminus L)$. That does not affect the value of the integral since the complement has measure zero, so:
\[
	r (g) = \lim_{p \to \infty} \frac{1}{p} \int_{X_n  (\Sigma \setminus L)} \psi (\gamma (g^p, \mathbf{x})) \dd \text{vol}.
\]
In the construction of $\gamma (g, \mathbf{x})$ we replace the geodesic segments from $z_j$ to $x_j$ and segments connecting $g (x_j)$ to $z_j$ by short paths in $\Sigma \setminus L$ with the same endpoints. Denote the result by $\gamma' (g, \mathbf{x}) \in MCG (\Sigma_n)$. 
For appropriate choice of connecting paths the distance
\[
   d\left(\gamma (g, \mathbf{x}), \gamma' (g, \mathbf{x})\right)
\]
is uniformly bounded in $g$ and $\mathbf{x}$ for any word metric (we provide a proof for this technical statement at the end of this section), 
hence $| \psi (\gamma (g, \mathbf{x})) - \psi (\gamma' (g, \mathbf{x}))|$ is bounded uniformly as well.
Therefore the difference will disappear under stabilization of $r$:
\[
	r (g) = \lim_{p \to \infty} \frac{1}{p} \int_{X_n  (\Sigma \setminus L)} \psi (\gamma (g^p, \mathbf{x})) \dd \text{vol} = 
				\lim_{p \to \infty} \frac{1}{p} \int_{X_n  (\Sigma \setminus L)} \psi (\gamma' (g^p, \mathbf{x})) \dd \text{vol}.
\]

Finally, both $g^p$ and the pushes of $\gamma' (g^p, \mathbf{x})$ preserve $L$, thus the essential curve $L$ is preserved under the composition.
That is, $\gamma' (g^p, \mathbf{x})$ is reducible. Bestvina-Fujiwara quasimorphisms vanish on reducible elements, hence the expression inside the integral is zero and $r (g) = 0$.

If $\Sigma \setminus L$ is not connected, there are few technical difficulties to overcome. First, a diffeomorphism $g$ may 
permute the connected components. But in this case we note that for an appropriate $l > 1$, $g^l$ brings all the components back, and since
$r (g) = \frac{r(g^l)}{l}$, it is enough to prove that $r(g^l) = 0$. In what follows we will assume that the connected components are not permuted by 
diffeomorphisms in question.
\begin{rem}
 In fact, if the essential curve $L$ is connected, a Hamiltonian diffeomorphism which preserves $L$ cannot move the connected components of $\Sigma \setminus L$. We analyze this case to allow curves $L$ with several connected components. We also wish to keep the argument valid for groups $\Symp_0$ and $\Symp$ rather than just $\Ham$. 
\end{rem}
$X_n (\Sigma \setminus L)$ has $k^n$ connected components (according to the location of each of the $n$ elements of $\mathbf{x} \in X_n (\Sigma \setminus L)$, $k$ is the number of connected components in $\Sigma \setminus L$). Denote these components by $X^1, \ldots, X^{k^n}$. As before,
\[
	r (g) = \lim_{p \to \infty} \frac{1}{p} \int_{X_n  (\Sigma \setminus L)} \psi (\gamma (g^p, \mathbf{x})) \dd \text{vol} = 
			   \lim_{p \to \infty} \sum_{i=1}^{k^n} \frac{1}{p} \int_{X^i} \psi (\gamma (g^p, \mathbf{x})) \dd \text{vol}.
\]

Fix $k^n$ $n$-tuples of marked points $\{\mathbf{z}^i \in X^i\}_{i=1}^{k^n}$ by picking one in each connected component.
We modify the construction of $\gamma (g, \mathbf{x})$ by replacing the basepoint $\mathbf{z}$ with $\mathbf{z}^i$ from the same connected component of 
$X_n (\Sigma \setminus L)$ as $\mathbf{x}$. 
Furthermore, we select short paths to and from $\mathbf{z}^i$ that avoid $L$. Denote the result as $\gamma' (g, \mathbf{x})$.
As explained a bit later, both changing the punctures and modification of paths result in a bounded deformation of $\gamma$:
$d\left(\gamma (g, \mathbf{x}), \gamma' (g, \mathbf{x})\right)$ and $| \psi (\gamma (g, \mathbf{x})) - \psi (\gamma' (g, \mathbf{x}))|$ 
are uniformly bounded, hence the effect of this deformation disappears under stabilization of the integral. That is,
\[
	r (g) = \lim_{p \to \infty} \sum_{i=1}^{k^n} \frac{1}{p} \int_{X^i} \psi (\gamma (g^p, \mathbf{x})) \dd \text{vol} = 
				\lim_{p \to \infty} \sum_{i=1}^{k^n} \frac{1}{p} \int_{X^i} \psi (\gamma' (g^p, \mathbf{x})) \dd \text{vol}.
\]
Similarly to the argument above, $\gamma' (g, \mathbf{x})$ preserves $L$, hence is reducible and $\psi (\gamma' (g, \mathbf{x})) = 0$.

\medskip 

It is left to analyze the effect of changing the paths and/or marked points in the construction of $\gamma (g, \mathbf{x})$.

We start with marked points: let $\mathbf{z}$ and $\mathbf{z}'$ be two $n$-tuples of punctures. Given 
$\mathbf{x} \in X_n (\Sigma)$ we follow the usual procedure to obtain $\gamma (g, \mathbf{x})$ and $\gamma' (g, \mathbf{x})$ - elements
of mapping class groups of $\Sigma$ with punctures at $\mathbf{z}$ and $\mathbf{z}'$, respectively.
In order to compare $\gamma (g, \mathbf{x})$ with $\gamma' (g, \mathbf{x})$ we need to identify the two mapping class groups.
This can be done by pushing points of $\mathbf{z}$ to $\mathbf{z}'$ along short geodesic segments, applying $\gamma' (g, \mathbf{x})$
and pushing the points $\mathbf{z}'$ back to $\mathbf{z}$ along short geodesics.
While this identification is not canonical (depends on the Riemannian metric), different identifications are conjugate hence do not affect
the values of the homogeneous quasimorphism $\psi$.

Under this identification, the pushes in the construction of $\gamma'$ go along broken geodesics from $\mathbf{z}$ to $\mathbf{z}'$
and then to $\mathbf{x}$, while $\gamma$ uses a straight push from $\mathbf{z}$ to $\mathbf{x}$. Same in the opposite direction when pushing
$g(\mathbf{x})$ to $\mathbf{z}$. That is, $\gamma'$ differs from $\gamma$ by a push along a $n$-tuple of short geodesic triangles
$\mathbf{z} \to \mathbf{z}' \to \mathbf{x} \to \mathbf{z}$ on the right and a similar $n$-tuple of triangles on the left. 
We claim that these triangular pushes have uniformly bounded norm in any metric (in particular, in any word metric).

Indeed, we push along a braid with $n$ strands where each strand traverses a geodesic triangle with short sides. Two short geodesic segments
with distinct endpoints can intersect at most once, therefore the image of this braid in $\Sigma$ has less than $9n^2$ intersection points.
As the result, the braid can be presented as a composition of elementary braids (those that swap two marked points or those that 
rotate one marked point along a simple non-contractible loop in the punctured surface, leaving the remaining strands fixed). 
The number of ingredients in such a decomposition is at most $9n^2 + n$. If one restricts the lengths of a noncontractible loop and 
of a trajectory used for swapping (each geodesic triangle consists of three short segments, hence its perimeter is bounded), 
the set of elementary braids which may show up in this decomposition is finite and independent of $g$ or $\mathbf{x}$. 
Hence the norm of this triangular push is at most $9n^2 + n$ times the maximal norm of a push along an elementary braid from the finite generating set.

At last, let $\gamma' (g, \mathbf{x})$ be a modification of $\gamma (g, \mathbf{x})$ where we keep the marked points $\mathbf{z}$ 
in place but replace the connecting paths by those that avoid the invariant curve $L$. We assume that for all $i$,
$x_i$ and $z_i$ belong to the same connected component of $\Sigma \setminus L$.
Let $c_i$ be a simple path connecting $z_i$ to $x_i$ and whose length is at most the intrinsic diameter of $\Sigma \setminus L$.
$c_i$ can be approximated in $\Sigma \setminus L$ by a piecewise geodesic path with short segments. Moreover, by standard compactness arguments, 
the number of segments in this construction can be uniformly bounded.
We use these piecewise geodesic paths to connect $z_i$ to $x_i$ and $g (x_i)$ back to $z_i$ in the construction of $\gamma'$.
The rest of the argument is similar to the previous case.
$\gamma'$ differs from $\gamma$ by a push of $\mathbf{z}$ along a piecewise geodesic braid on the left and on the right.
Geodesic segments are short and their total number is bounded, hence there is a uniform bound on the number of intersections and 
these braids can be presented as a bounded composition of elementary braids.

\section{The annulus}

We prove Theorem~\ref{t:annulus}.

\subsection{Tools}

Let $F_t : \Sigma \to \RR$, $t \in [0, 1]$ be a time-dependent smooth function with compact support in the interior of $\Sigma$. We define
$\widetilde{\Cal} (F_t) = \int_0^1 \left( \int_\Sigma F_t \omega \right) \dd t$. If the symplectic form $\omega$ is exact 
(this is the case for an annulus or a disk),
$\widetilde{\Cal}$ descends to a homomorphism $\Cal_\Sigma: \Ham(\Sigma) \to \RR$ which is called 
the Calabi homomorphism.

\medskip

Let $\Sigma$ be a symplectic surface of genus zero.
Given a compactly supported smooth function $F: \Sigma \to \RR$, the \emph{Reeb graph} $T_F$ is defined as the set of connected
components of level sets of $F$ (for a more detailed description we refer the reader to ~\cite{En-Po:calqm}). 
For a generic Morse function $F$ (saying `Morse', we mean that the restriction of $F$ to the interior of its support 
is a Morse function) this set, equipped with topology 
induced by the projection $\pi_F: \Sigma \to T_F$, is homeomorphic to 
a tri-valent tree. We endow $T_F$ with a positive measure given by $\mu (X) = \int_{\pi_F^{-1}(X)} \omega$ 
for all $X \subseteq T_F$ with measurable $\pi_F^{-1}(X)$. 
In the case of the annulus $\Sigma = S^1 \times [0, 1]$, $\pi_F (S^1 \times \{0\})$ will be 
referred to as the \emph{bottom root} of $T_F$ and $\pi_F (S^1 \times \{1\})$ as the \emph{top root}. The shortest path connecting the roots of $T_F$ will be called a \emph{stem}.

A point $x_{m} \in T_F$ is a \emph{median} of $T_F$ if all connected components of $T_F \setminus \{x_m\}$ have measure at most $\frac{\Area (\Sigma)}{2}$.
A median always exists and is unique (see ~\cite{En-Po:calqm}). The set $\pi_F^{-1} (x_m)$ will be called the \emph{median} with respect to $F$.
Suppose $\Sigma = S^1 \times [0, 1]$, we define \emph{percentile} sets in analogy to the median.
Let $h \in [0, 1]$. $x_h \in T_F$ is an \emph{$h$-percentile} of $T_F$ if the top and the bottom roots belong to different connected components of $T_F \setminus \{x_h\}$ 
and the connected component of the bottom root has measure $h \cdot \Area (\Sigma)$. 
The set $\pi_F^{-1} (x_h)$ is an $h$-percentile with respect to $F$. 

Clearly, percentiles correspond to points $x$ in the stem of $T_F$ and the percentile value increases monotonically 
along the stem. Unlike the median, if $T_F$ is not homeomorphic to an interval (that is, has `branches' besides the stem), 
$h$-percentiles do not exist for certain $h \in [0, 1]$. Each branch corresponds to a `gap' (missing interval) in the set of percentile values. Length of the gap is given by the 
measure of the branch normalized by $\Area (\Sigma)$. If an $h$-percentile exists, it is unique. 
The $\frac{1}{2}$-percentile (if it exists) coincides with the median. For a generic $F$ this corresponds to the case when the median set of $F$ is a non-contractible circle.
Using a standard Morse-theoretic argument, we conclude with the following observation: percentile sets are not contractible in $S^1 \times [0, 1]$. 
The set $A_F$ of points that are not percentiles of $T_F$
is the union of branches that grow out of the stem of $T_F$. The set $\pi^{-1} (A_F)$ is the union of topological disks corresponding to these branches.

\medskip

In ~\cite{En-Po:calqm} the authors describe construction of a homogeneous quasimorphism $$\Cal_{S^2} : \Ham (S^2) \to \RR.$$
It has the following properties: $\Cal_{S^2}$ is Hofer-Lipschitz $$|\Cal_{S^2}(\phi)| \leq \Area (S^2) \cdot \| \phi \|.$$
In the case when $\phi \in \Ham(S^2)$ is supported in a disk $D$ which is displaceable 
in $S^2$, $\Cal_{S^2} (\phi) = Cal_D (\phi \big|_D)$.
Moreover, for a $\phi \in \Ham(S^2)$ generated by an autonomous function $F: S^2 \to \RR$, 
$\Cal_{S^2} (\phi)$ can be computed in the following way. 
Let $x$ be the median of $T_F$ and $X = \pi_F^{-1} (x)$ be the corresponding subset of $S^2$.
Then 
\[
	\Cal_{S^2}(\phi) = \int_{S^2} F \omega - \Area (S^2) \cdot F(X).
\]

Let $\Sigma = S^1 \times [0, 1]$ be an annulus equipped with the standard symplectic form so that $\Area (\Sigma) = 1$.
We embed $\Sigma$ into a sphere $S^2_{a,b}$ of area $1 + a + b$ by gluing a disk of area $a$ to $S^1 \times \{0\}$ and a disk of area $b$ to $S^1 \times \{1\}$.
Denote this embedding by $i_{a, b} : \Sigma \to S^2_{a,b}$.
Let $$r_{a,b} = \frac{1}{1+a+b} \cdot \left( \Cal_\Sigma - i_{a,b}^* \Cal_{S^2_{a,b}} \right)$$ be the normalized difference between the Calabi homomorphism on $\Sigma$ and the pullback of the Calabi quasimorphism of $S^2_{a,b}$.
Note that $r_{a,b}$ vanishes on Hamiltonians $g$ supported in a disk $D \subset \Sigma$ of area $\frac{1+a+b}{2}$. Indeed, $i_{a,b} (D)$ is displaceable in $S^2_{a,b}$ thus 
$$\Cal_{S^2_{a,b}} (i_{a,b,*} g) = \Cal_D \left(g \big|_D \right) = \Cal_\Sigma (g).$$ This implies that $r_{a,b}$ is continuous in the $C^0$-topology (see ~\cite{En-Po-Py:qm-continuity}).

Let $F: \Sigma \to \RR$ be a Hamiltonian function, $f$ its time-$1$ map and suppose that $-1 \leq b-a \leq 1$ or, equivalently, $h: = \frac{1+b-a}{2} \in [0, 1]$. If $F$ admits  
the $h$-percentile set $X_h$, it is mapped by $i_{a, b}$ to the median set of $i_{a, b, *} F : S^2_{a,b} \to \RR$, therefore 
$r_{a,b} (f) = F(X_h)$. This makes the quasimorphisms $r_{a,b}$ a useful tool to extract information about the Reeb graph of a 
Hamiltonian function.

\subsection{Construction}
We construct a non-autonomous Hamiltonian on $\Sigma$. Later we will show that it induces a non-autonomous deformation
on $L$. Let $F : S^1 \times [0, 1] \to \RR$ be a Hamiltonian function
given by $F (\theta, s) = s$ when $s \in [0.01, 0.99]$ and extended to the rest of $\Sigma$ in arbitrary way. The time-$t$ map $f_t$ of $F$ rotates the annulus $A = S^1 \times [0.01, 0.99]$
by $t$ in the $S^1$ coordinate. Let $D \subset A$ be a disk of area $0.8$ and $\Phi : \Sigma \to \RR$ be a smooth function which equals $1$ in $D$ and is supported in a disk of area $0.9$ inside $A$. (That is, $\Phi$ is a smooth approximation of the indicator function of $D$.)
The time-$t$ map $\phi_t$ fixes $D$ pointwise but the flow induces a fast rotation outside $\partial D$.
Pick large independent parameters $T, \tau \in \NN$ and consider $g_{T, \tau} := f_T \circ \phi_\tau$. Assuming $T$ is an integer, $f_T$ translates 
the subannulus $A$ precisely $T$ times around $S^1$, hence fixes $A$ pointwise.
$\phi_\tau$ is supported in $A$, hence $f_T$ and $\phi_\tau$ commute.

\begin{figure}[!htbp]
\begin{center}
\includegraphics[width=0.7\textwidth]{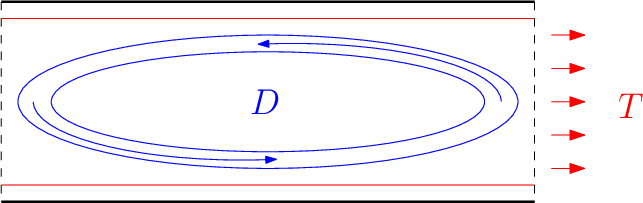}
\caption{$g_{T, \tau}$}
\end{center}
\end{figure}

We claim that $g_{T, \tau}$ is not autonomous. Assume by contradiction that it is generated as the time-$1$ map of a Hamiltonian function $H : \Sigma \to \RR$. 
Suppose first that $H$ is generic, that is, $H$ admits a Reeb tree $T_H$.
We compute the values of $H$ at its percentile sets in two different ways: first, pick $h \in [0.01, 0.99]$. Let $a = 1$ and $b = 2h$ which satisfy $h = \frac{1+b-a}{2}$.
\[
  r_{a,b} (g_{T, \tau}) = r_{a,b} (f_T) + r_{a,b} (\phi_\tau) = h T.
\]
The first equality holds because $f_T$ and $\phi_\tau$ commute. $r_{a,b} (f_T) = h T$ since $Y_h = S^1 \times \{h\}$ is the $h$-percentile for $F$ and $F(Y_h) = h$.
$r_{a,b} (\phi_\tau) = 0$ as the support of $\Phi$ becomes displaceable in $S^2_{a,b}$.
Therefore, if the $h$-percentile $X_h$ exists for $H$, $H (X_h) = h T$.

We perform another computation: fix $h' \in [0.2, 0.8]$. Let $a' = 0.8 - h'$ and $b' = h' - 0.2$. Once again, $h' = \frac{1+b'-a'}{2}$ and
\[
  r_{a',b'} (g_{T, \tau}) = r_{a',b'} (f_T) + r_{a',b'} (\phi_\tau) = h' T + \tau.
\]
$r_{a',b'} (f_T) = h' T$ as before but $r_{a',b'} (\phi_\tau) = \tau$ since $i_{a', b'}$ embeds $\Sigma$ into a sphere of area $1 + a' + b' = 1.6$. So the image of the disk $D$ becomes the median set
for $i_{a', b',*} \Phi$, thus $r_{a',b'} (\phi_\tau)$ can be computed explicitly.
The calculation shows that if the $h'$-percentile $X_{h'}$ exists for $H$, $H (X_{h'}) = h' T + \tau$.

This contradicts the previous result, hence $h$-percentiles do not exist for $h$ in the interval $[0.2, 0.8]$. 
That is, $T_H$ has one or several branches with total measure at least $0.6$. In fact, there must be a single branch of measure at least $0.6$: 
if there are several branches growing out of different points of the stem, there will be intermediate $h$-percentiles which correspond to stem points between the branches. 
In our situation it is not the case. If there are two branches or more growing from the same stem point (which is possible in a non-generic situation), we may perturb $H$ in 
the $C^\infty$-topology and separate the branches. Intermediate percentiles will appear after such perturbation. However, our quasimorphisms $r_{a,b}$ are $C^0$-continuous, so a small perturbation
will not resolve the discrepancy $\tau$ between the results of two computations. 

As a corollary, there must be a branch $B \subset T_H$ with measure at least $0.6$. $D_B = H^{-1} (B)$ is a topological disk in $\Sigma$ of area at least $0.6$ which is an invariant set 
for the flow of $H$. Intuitively, points in $D_B$ have rotation number $0$ with respect to the $S^1$ coordinate
(all points with non-zero rotation number are mapped to the stem). However, most points in $\Sigma$ (up to a subset of area $0.02$) have 
rotation number $T$ under $g_{T, \tau}$, which gives a contradiction. 

\medskip

We reproduce this contradiction using more powerful tools. In ~\cite{Kh:h-spec}, Theorem 2, the author constructs a quasimorphism $\rho_{0.6} : \Ham (\Sigma) \to \RR$ which is $C^0$-continuous and has the following property.
Suppose $g \in \Ham (\Sigma)$ has an invariant disk of area $0.6$ or more, then $\rho_{0.6} (g)$ computes the rotation number (along the $S^1$ coordinate) of points in this disk.
($\rho_{0.6}$ is constructed as a certain combination of Calabi quasimorphisms pulled back from $S^2$ similarly to the construction of $r_{a,b}$.)
Therefore, 
\[
	\rho_{0.6} (g_{T, \tau}) = \rho_{0.6} (f_T) + \rho_{0.6} (\phi_\tau) = T.
\]
$\rho_{0.6} (f_T) = T$ since $f_T$ rotates the annulus $A$ $T$ times around, the same is true for any disk of area $0.6$ in $A$.
$\rho_{0.6} (\phi_\tau) = 0$ since $D$ is a stationary disk of area $0.8$.

This shows that large invariant disks of $g_{T, \tau}$ (if they exist) have rotation number $T$.
On the other hand, invariant disks of an autonomous flow must have rotation number zero. Therefore no such disks exist, 
so the Hamiltonian function $H$ cannot have a large branch. This is a contradiction to the first part of the 
argument where we established existence of a branch $B$. Hence $g_{T, \tau}$ is not autonomous.

If $H$ which is supposed to generate $g_{T, \tau}$ is extremely non-generic and its Reeb graph does not exist, we may perturb it and argue as before, since 
the quasimorphisms $r_{a,b}$ and $\rho_{0.6}$ used as tools to arrive to a contradiction are $C^0$-continuous.

\begin{rem}
  $g_{T, \tau}$ is not autonomous in $\Ham(\Sigma)$ but is a composition of two autonomous maps.
  
  However, if one allows Hamiltonian flows and diffeomorphisms in $\Sigma$ whose support is not restricted to the interior, $g_{T, \tau}$ 
  becomes autonomous in this extended group and, in particular, has entropy zero.
  To see this, note that for an integer $T$, the map $f_T$ which rotates the inner subannulus $A$ around the $S^1$ coordinate, 
  can be generated by another autonomous flow $\tilde{f}_t$. Namely, the one that fixes $A$ pointwise and rotates a tubular neighborhood of $\partial \Sigma$
  in the opposite direction. This flow is generated by $\widetilde{F}(\theta, s) = F (\theta, s) - s$.
  Since $\phi_t$ is supported inside $A$, $\tilde{f}_t$ and $\phi_t$ have disjoint supports, commute and can be combined into a single autonomous flow which
  generates $g_{T, \tau} = f_T \circ \phi_\tau$.
\end{rem}
 
 The obstruction for $g_{T, \tau}$ to be autonomous consists of two ingredients:
  \begin{itemize}
	  \item
		 $\forall h \in [0.2, 0.8]. \, r_{a',b'} (g_{T, \tau}) -  r_{a,b} (g_{T, \tau}) = \tau \neq 0$, 
		 hence no percentiles exist in the interval $[0.2, 0.8]$. Therefore there must be a branch of area at least $0.6$.	  
	  \item
	   $\rho_{0.6} (g_{T, \tau}) = T \neq 0$.
	   Therefore $g_{T, \tau}$ cannot have a large invariant disk with rotation number zero, hence no large branches for 
	   the generating function of $g_{T, \tau}$.
  \end{itemize}
  The quasimorphisms $\rho$ and $r$ used in the argument are Hofer-Lipschitz, 
  hence this obstruction persists under deformations of $g_{T, \tau}$ whose Hofer norm is less than $\min (T, \tau)$ divided by appropriate Lipschitz constants.
  This provides a lower bound for the Hofer distance between $g_{T, \tau}$ and the set of autonomous Hamiltonians. In particular,
  $g_{T, \tau}$ arrives arbitrarily far away from autonomous diffeomorphisms if we let $T, \tau \to \infty$. 
  
\begin{rem}
  We compare $g_{T, \tau}$ with the egg-beater maps of Polterovich and Shelukhin (see ~\cite{Po-Sh:egg-beater}). An egg-beater map can also be constructed arbitrarily far away in Hofer's metric 
  from any autonomous Hamiltonian.  But it is constructed on surfaces of higher genus, it 
  is highly chaotic and has positive entropy, which is very different from our example. In addition, egg-beaters stay far away also from powers of Hamiltonian diffeomorphisms
  while $g_{N, N} = g_{1,1}^N$.
  
  On the other hand, $h$-percentiles and invariants computed by the quasimorphisms $r$ and $\rho$ can be expressed in terms of persistence modules, so our
  methods may have common background with those of ~\cite{Po-Sh:egg-beater}.
\end{rem}

\begin{rem}
  Another direction for comparison is quasimorphisms on surfaces that vanish on autonomous diffeomorphisms (see ~\cite{Bra-Ma:ent-quas} and a series of earlier works
  ~\cite{Bra-Ke:auton-disc, Bra-Ke-Sh:auton-torus, Bra:biinv-quas}). Both approaches use quasimorphisms as tools. However, the quasimorphisms used here do not vanish on autonomous
  Hamiltonians, hence cannot be used directly to prove the desired result or to construct Hamiltonians that are far from the identity in the autonomous norm. On the positive side,
  our quasimorphisms are Hofer-Lipschitz and descend to invariants of curves in $S^1 \times [0, 1]$ (which is not the case in \cite{Bra-Ma:ent-quas}).
\end{rem}

\medskip

We are ready to prove Theorem~\ref{t:annulus}.
Let $L = \{0\} \times [0, 1] \cup \{\frac{1}{2}\} \times [0, 1] \subset \Sigma$, $L' = g_{T, \tau}(L)$. 
We show that quasimorphisms $r_{a,b}$ and $\rho_{0.6}$ descend as ill-defined invariants to the orbit $O_L$.

Let $g$ be a Hamiltonian in the stabilizer $S$, that is, $g L = L$. We may perturb $g$ by a Hamiltonian $h$ supported in a small neighborhood of $L$ so that $hg$ fixes a neighborhood of 
$L$ pointwise. $hg = g_1 \circ g_2$ splits into a composition of two Hamiltonian diffeomorphisms: $g_1$ supported in $(0, \frac{1}{2}) \times [0, 1]$ and $g_2$ in $(\frac{1}{2}, 1) \times [0, 1]$.
Both supported in a topological disk of area $\frac{1}{2} \leq \frac{1 + a + b}{2}$, hence $r_{a,b} (g_1) = r_{a,b} (g_2) = 0$. $g_1, g_2$ commute (their supports are disjoint), $r_{a,b}$ is homogeneous,
hence $$r_{a,b} (hg) = r_{a,b}(g_1) + r_{a,b} (g_2) = 0.$$ $r_{a,b} (h) = 0$ by the same reason, which implies $|r_{a,b} (g)| \leq D_{r_{a,b}}$. That is, the restriction of $r_{a,b}$ to
the subgroup $S$ is bounded. $r_{a,b} \big|_S$ is a homogeneous bounded quasimorphism, hence it is identically zero. 

Similarly, $hg$ fixes a large topological disk given by removing a neighborhood of $\{0\} \times [0, 1]$ from $\Sigma$. It has rotation number zero, hence $\rho_{0.6} (hg) = 0$.
The same is true for $h$, so $\rho_{0.6} (h) = 0$. We continue as before: $|\rho_{0.6} (g)| \leq D_{\rho_{0.6}}$ and the quasimorphism vanishes on $S$.

Therefore all estimates and computations of quasimorphisms carried out for $g_{T, \tau}$ remain valid for 
the equivalence class $[g_{T, \tau}] = g_{T, \tau} S \in \Ham/S$ up to a compensation of ambiguity (which is bounded by the defects).
In another formulation, they descend to $g_{T, \tau}(L) \in O_L \simeq \Ham/S$.
Indeed, let $g' \in \Ham (\Sigma)$ such that $g' (L) = g_{T, \tau} (L)$. 
Given $h \in [0.2, 0.8]$ pick $a, b, a', b'$ adjusted to $h$ as in the beginning of the section.
$g'$ differs from $g_{T, \tau}$ by an element of $S$, hence 
\[
   r_{a',b'} (g') - r_{a,b} (g') > r_{a',b'} (g_{T, \tau}) - r_{a,b} (g_{T, \tau}) - D_{r_{a',b'}} - D_{r_{a,b}} =  \tau - D_{r_{a',b'}} - D_{r_{a,b}} .
\]
As before, we deduce that for $\tau$ large enough the autonomous function which generates $g'$ (if it exists) must have a large branch $B$. But
\[
	\rho_{0.6} (g') > \rho_{0.6} (g_{T, \tau}) - D_{\rho_{0.6}} =  T - D_{\rho_{0.6}} .
\]
If there is a large branch for $g'$, $\rho_{0.6} (g')$ will compute its rotation number which must be zero 
(all branches are stationary under the flow). This is a contradiction. That is, our obstruction for autonomous Hamiltonians
applies to all $\{g' \in \Ham (\Sigma) \, \big| \, g' (L) = L' \}$.

Due to Hofer-Lipschitz property of the quasimorphisms, this obstruction persists under deformations of $L' = g_{T, \tau}(L)$ unless
the deformation has Hofer's norm comparable to $\min (T, \tau)$. The last part of the theorem holds if we pick $T, \tau >> R$.

\bibliography{bibliography}

\begin{thebibliography}{{Cal}}

\bibitem[BF]{Be-Fu:bounded-cohomology}
Mladen Bestvina and Koji Fujiwara.
\newblock Bounded cohomology of subgroups of mapping class groups.
\newblock {\em Geom. Topol.}, 6(1):69--89, 2002.

\bibitem[BK]{Bra-Ke:auton-disc}
Michael Brandenbursky and Jarek K\k{e}dra.
\newblock On the autonomous metric on the group of area-preserving
  diffeomorphisms of the 2-disc.
\newblock {\em Algebraic \& Geometric Topology}, 13, 07 2012.

\bibitem[BKS]{Bra-Ke-Sh:auton-torus}
Michael Brandenbursky, Jarek K\k{e}dra, and Egor Shelukhin.
\newblock On the autonomous norm on the group of {H}amiltonian diffeomorphisms
  of the torus.
\newblock {\em Communications in Contemporary Mathematics}, 20(02):1750042,
  2018.

\bibitem[BM]{Bra-Ma:ent-quas}
Michael Brandenbursky and Micha\l{} Marcinkowski.
\newblock Entropy and quasimorphisms.
\newblock {\em Journal of Modern Dynamics}, 15:143--163, 2019.

\bibitem[Bra]{Bra:biinv-quas}
M.~Brandenbursky.
\newblock Bi-invariant metrics and quasi-morphisms on groups of {H}amiltonian
  diffeomorphisms of surfaces.
\newblock {\em Int. J. of Math.}, 26(9), 2013.

\bibitem[{Cal}]{Ca:scl}
Danny {Calegari}.
\newblock {\em {scl.}}, volume~20.
\newblock Tokyo: Mathematical Society of Japan, 2009.

\bibitem[EP]{En-Po:calqm}
M.~Entov and L.~Polterovich.
\newblock Calabi quasimorphism and quantum homology.
\newblock {\em Int. Math. Res. Not.}, 2003(30):1635--1676, 2003.

\bibitem[EPP]{En-Po-Py:qm-continuity}
M.~Entov, L.~Polterovich, and P.~Py.
\newblock On continuity of quasimorphisms for symplectic maps.
\newblock In {\em Perspectives in Analysis, Geometry, and Topology (a volume
  dedicated to Oleg Viro's 60th birthday)}, volume 296 of {\em Progress in
  Mathematics}. Birkh\"{a}user, Basel, 2012.

\bibitem[Kha]{Kh:h-spec}
M.~Khanevsky.
\newblock Hofer's length spectrum of symplectic surfaces.
\newblock {\em J. of Modern Dynamics}, 9(1):219--235, 2015.

\bibitem[PS]{Po-Sh:egg-beater}
Leonid Polterovich and Egor Shelukhin.
\newblock Autonomous {H}amiltonian flows, {H}ofer's geometry and persistence
  modules.
\newblock {\em Selecta Mathematica}, 22(1):227--296, Jan 2016.

\bibitem[SRS]{RSS:combi-floer-homology}
Vin Silva, Joel Robbin, and Dietmar Salamon.
\newblock Combinatorial {F}loer homology.
\newblock {\em Memoirs of the American Mathematical Society}, 230, 05 2012.

\end{thebibliography}

\end{document}